\documentclass{article}
\usepackage{enumitem}
\usepackage{latexsym,amsfonts,amssymb}
\usepackage[latin1]{inputenc}
\usepackage{amsmath}

\title{\sc Normal subgroups whose conjugacy class graph has diameter three.}

\author{Antonio Beltr\'an\\
\footnotesize
Departamento de Matem\'aticas,\\
\footnotesize Universidad Jaume I, \footnotesize
12071 Castell\'on, Spain\\
\footnotesize
e-mail: abeltran@mat.uji.es\\
\\
Mar\'{\i}a Jos\'e Felipe\\
\footnotesize
Instituto Universitario de Matem\'atica Pura y Aplicada,\\
\footnotesize Universidad Polit\'ecnica de Valencia, \footnotesize
46022 Valencia, Spain\\
\footnotesize e-mail: mfelipe@mat.upv.es \\
\\ Carmen Melchor\\
\footnotesize
Departamento de Educaci\'on,\\
\footnotesize Universidad Jaume I, \footnotesize
12071 Castell\'on, Spain\\
\footnotesize
e-mail: cmelchor@uji.es      }

\date{}
\begin{document} \maketitle

\begin{abstract}
Let $G$ be a finite group and $N$ a normal subgroup of $G$. We determine the structure of $N$ when the diameter of the graph associated to the $G$-conjugacy classes contained in $N$ is as large as possible, that is, is equal to three.



\end{abstract}

\bigskip

\section{Introduction}
Let $G$ be a finite group and let $N$ be a normal subgroup of $G$ and let $x\in N$. We denote by $x^{G}=\lbrace x^{g} \mid g \in G \rbrace$ the $G$-conjugacy class of $x$. Let $\Gamma_{G}(N)$ be the graph associated to these $G$-conjugacy classes, which was defined in \cite{Nuestro1} as follows: its vertices are the $G$-conjugacy classes of $N$ of cardinality bigger than 1, that is, $G$-classes of elements lying in $N\setminus ({\bf Z}(G)\cap N)$, and two of them are joined by an edge if their sizes are not coprime. It was proved in \cite{Nuestro1} that $d(\Gamma_{G}(N))\leq 3$, where $d(\Gamma_{G}(N))$ denotes the diameter of the graph. In this note we analyze the structure properties of $N$ when $d(\Gamma_{G}(N))=3 $.\\

The above graph extends the ordinary graph, $\Gamma(G)$, which was formerly defined in \cite{BerHerMann}, and whose vertices are the non-central conjugacy classes of $G$ and two vertices are joined by an edge if their sizes are not coprime. The graph $\Gamma_G(N)$ can be viewed as the subgraph of $\Gamma(G)$ induced by those vertices of $\Gamma(G)$ which are vertices in $\Gamma_G(N)$. This fact does not allow us however, to obtain directly properties of the graph of $G$-classes.\\

Concerning ordinary classes, L.S. Kazarin characterizes in \cite{Kazarin} the structure of a group $G$ having two ``isolated classes". We recall that a group $G$ is said to have isolated classes if there exist elements $x, y \in G$
such that every element of $G$ has a conjugacy class size coprime to either $|x^{G}|$ or $|y^{G}|$. Particularly, Kazarin determined the structure of those groups $G$ with $d(\Gamma(G))=3$. On the other hand, the disconnected graph was studied by Bertram, Herzog and Mann in \cite{BerHerMann}. It should be noted that similar results have also been studied for other graphs. In \cite{Dolfi}, Dolfi defines the graph $\Gamma'(G)$ whose vertices are the elements of the set of all primes which occur as divisors of the lengths of the conjugacy classes of $G$, and two vertices $p,q$ are joined by an edge if there exists a conjugacy class in $G$ whose length is a multiple of $pq$. In \cite{DolfiCasolo} Dolfi and Casolo describe all finite groups $G$ for which $\Gamma'(G)$ is connected and has diameter three. \\

We remark that the primes dividing the $G$-conjugacy class sizes do not need to divide $|N|$. This specially occurs when $N$ is abelian and it is non-central in $G$ and consequently, we may have no control on this set of primes. For this reason, we observe that new cases appear when dealing with $G$-classes which are not contemplated in the ordinary case. The main result of this note is Theorem A. From now on, if $G$ is a finite group we denote by $\pi(G)$ the set of primes dividing $|G|$.\\

{\bf Theorem A.} {\it  Let $G$ be a finite group and $N\unlhd G$. Suppose that $x^G$ and $y^G$ are two non-central $G$-conjugacy classes of $N$ such that any $G$-conjugacy of $N$ has size  coprime with $|x^G|$ or $|y^G|$. Let $\pi_{x}=\pi(|x^G|)$, $\pi_{y}=\pi(|y^G|)$ and $\pi=\pi_{x}  \cup \pi_{y}$. Then, $N=$ {\rm \textbf{O}}$_{\pi'}(N)\times ${\rm \textbf{O}}$_{\pi}(N)$ with $x, y \in$ {\rm \textbf{O}}$_{\pi}(N)$, which is either a quasi-Frobenius group with abelian kernel and complement or {\rm \textbf{O}}$_{\pi}(N)=P\times A$ with $A\leq$ {\rm{\textbf{Z}}}$(N)$ and $P$ a $p$-group for a prime $p$.}\\

Notice that in the conditions of Theorem A, we have two possibilities:  $d(\Gamma_{G}(N))\leq 2$ or  $d(\Gamma_{G}(N))=3$. In the former case, we certainly have that the graph is disconnected, and then the structure of $N$ is already determined by Theorem E of \cite{Nuestro1}. However, we slightly improve this result in Corollary B. In the second case, the graph is connected since when the graph $\Gamma_{G}(N)$ is disconnected, it was proved that each connected component is a complete graph (Theorem B of \cite{Nuestro1}). Therefore, we deduce the following consequences for each of these cases.\\

{\bf Corollary B.} {\it Let $G$ be a finite group and $N\unlhd G$. Suppose that $\Gamma_{G}(N)$ is disconnected and let $x,y \in N$ such that $(|x^{G}|, |y^{G}|)=1$. Set $\pi=\pi(|x^{G}|) \cup \pi(|y^{G}|)$. Then $x$, $y\in $ {\rm \textbf{O}}$_{\pi}(N)$, $N= $ {\rm \textbf{O}}$_{\pi'}(N)\times ${\rm \textbf{O}}$_{\pi}(N)$ with {\rm \textbf{O}}$_{\pi'}(N)\subseteq$ {\rm \textbf{Z}}$(G)$, and either {\rm \textbf{O}}$_{\pi}(N)$ is a quasi-Frobenius group with abelian kernel and complement or {\rm \textbf{O}}$_{\pi}(N)=P\times A$ with $A\leq$ {\rm{\textbf{Z}}}$(G)$ and $P$ a $p$-group for a prime $p$.}\\

{\bf Corollary C.} {\it Let $G$ be a finite group and $N\unlhd G$. Suppose that $\Gamma_{G}(N)$ is connected with $d(\Gamma_{G}(N))=3$. Let us consider $x,y \in N$ such that $d(x^{G}, y^{G})=3$. Set $\pi=\pi(|x^{G}|) \cup \pi(|y^{G}|)$. Then $x$, $y\in $ {\rm \textbf{O}}$_{\pi}(N)$, $N= $ {\rm \textbf{O}}$_{\pi'}(N)\times ${\rm \textbf{O}}$_{\pi}(N)$ where  either {\rm \textbf{O}}$_{\pi}(N)$ is a quasi-Frobenius group with abelian kernel and complement or {\rm \textbf{O}}$_{\pi}(N)=P\times A$ with $A\leq$ {\rm{\textbf{Z}}}$(N)$ and $P$ a $p$-group for a prime $p$.}\\

\section{Proofs}

First, we state three elementary results which are need to prove the main result.\\

{\bf Lemma 1.} {\it Let $G$ a $\pi$-separable group. Then the conjugacy class length of every $\pi$-element of $G$ is a $\pi$-number if and only if $G=H\times K$, where $H$ and $K$ are a Hall $\pi$-subgroup and a $\pi$-complement of $G$, respectively.}\\

{\it Proof.} This is Lemma 8 of \cite{piseparable}. $\Box$ \\

In the particular case in which $\pi=p'$, the complement of some prime $p$, Lemma 1 is true without assuming $p$-separability (or equivalently $p$-solvability).\\

{\bf Lemma 2.} {\it Every $p'$-element of a group $G$ has index prime to $p$, for some prime $p$, then the Sylow $p$-subgroup of $G$ is a direct factor of $G$.}\\

{\it Proof.} This is Lemma 1 of \cite{Camina}. $\Box$\\

{\bf Lemma 3.} {\it Let $G$ be a finite group and $N\unlhd G$. Let $B=b^{G}$ and $C=c^{G}$ be two non-central $G$-conjugacy classes of $N$. If $(|B|,|C|)=1$. Then
\begin{enumerate}[label={\rm \alph*.}]
\item {\rm \textbf{C}}$_{G}(b)${\rm \textbf{C}}$_{G}(c)=G$.
\item $BC=CB$ is a non-central $G$-class of $N$ and $|BC|$ divides $|B| |C|$.
\item Suppose that $d(B,C)\geq3$ and $|B|<|C|$. Then $|BC|=|C|$ and $CBB^{-1}=C$. Furthermore, $C\langle BB^{-1}\rangle=C$, $\langle BB^{-1}\rangle\subseteq\langle CC^{-1}\rangle$ and $|\langle BB^{-1}\rangle |$ divides $|C|$.
\end{enumerate}}

{\it Proof.} This is Lemma 2.1 of \cite{Nuestro1}.$\Box$\\

{\it Proof of Theorem A.} We proceed by induction on $|N|$. Notice that the hypotheses are inherited by every normal subgroup in $G$ which is contained in $N$ and contains $x$ and $y$. By using the primary decomposition we can assume that $x$ and $y$ have order a power of two primes, say $p$ and $q$, respectively.\\

{\it Step 1. $q=p$ if and only if $xy=yx$.}\\

Suppose that $xy=yx$ and that $p\neq q$. Observe that ${\rm\textbf{C}}_{G}(xy)={\rm\textbf{C}}_{G}(x)\cap {\rm\textbf{C}}_{G}(y)$ and consequently, $|x^{G}|$ divides $|(xy)^{G}|$ and $|y^{G}|$ divides $|(xy)^{G}|$. Thus, we obtain a $G$-conjugacy class connected with $x^{G}$ and $y^{G}$, which is a contradiction with the hypotheses. Conversely, suppose that $p=q$. We know that $p$ cannot divide either $|x^{G}|$ or $|y^{G}|$. Furthermore, the hypotheses imply that $(|x^G|, |y^G|)=1$, so we have $G={\rm\textbf{C}}_{G}(x)\textbf{C}_{G}(y)$ and $|x^{G}|=|G:{\rm \textbf{C}}_{G}(x)|=|{\rm\textbf{C}}_{G}(y): {\rm\textbf{C}}_{G}(x)\cap {\rm\textbf{C}}_{G}(y)|$. Now, since $y$ is a $p$-element in ${\rm \textbf{Z}}({\rm\textbf{C}}_{G}(y))$, we deduce that $y \in {\rm\textbf{C}}_{G}(x) \cap {\rm\textbf{C}}_{G}(y)$ and hence $xy=yx$.\\

{\it Step 2. $p,q\in \pi$.}\\

 We define $K={\rm\textbf{C}}_{G}(x)\cap {\rm\textbf{C}}_{G}(y)$. First, we assume that $p\neq q$ and $xy\neq yx$. We have $|G:K|=|G:{\rm \textbf{C}}_{G}(x)||{\rm \textbf{C}}_{G}(x):{\rm \textbf{C}}_{G}(x)\cap {\rm \textbf{C}}_{G}(y)|=|x^{G}||y^{G}|$, which is a $\pi$-number. Since $x\in {\rm\textbf{Z}}({\rm\textbf{C}}_{G}(x))$ and $x$ is a $p$-element but $x\not\in K$, we know that $p$ divides $|{\rm\textbf{C}}_{G}(x):K|=|y^{G}|$. This means that $p\in \pi_{y}$. Similarly we obtain that $q$ divides $|x^{G}|$, that is, $q\in \pi_{x}$. Consequently, $p,q \in \pi$.\\

Suppose now that $p=q$ and $xy=yx$. Let us see that $p\in \pi$. We denote $X=x^{G}$ and $Y=y^{G}$ and we assume for instance that $|X|>|Y|$. By hypothesis, the distance between $X$ and $Y$ in $\Gamma_{G}(N)$ is bigger than or equal to 3. We can apply Lemma 3(c) and we obtain $X\langle YY^{-1}\rangle=X$, $\langle YY^{-1}\rangle\subseteq\langle XX^{-1}\rangle$ and $|\langle YY^{-1}\rangle |$ divides $|X|$. On the other hand, since $G={\rm \textbf{C}}_{G}(x){\rm \textbf{C}}_{G}(y)$ we have $X\subseteq {\rm \textbf{C}}_{G}(y)$. As a result, $\langle YY^{-1}\rangle\subseteq \langle XX^{-1}\rangle\subseteq {\rm \textbf{C}}_{G}(y)$. In particular, if we take $z=y^{g}\neq y$, for some $g\in G$, we have $w=zy^{-1}\in \langle YY^{-1}\rangle \subseteq {\rm \textbf{C}}_{G}(y)$, so $[z,y]=1$. We obtain that $w$ is a non-trivial $p$-element and, since $p$ divides $|\langle YY^{-1}\rangle|$, which divides $|X|$, we conclude that $p\in \pi_{x}$. If $|Y|>|X|$ we can argue similarly to get $p\in \pi_{y}$.\\

{\it Step 3. We can assume that $N/${\rm \textbf{Z}}$(N)$ is neither a $p$-group nor a $q$-group (particularly, we can assume that $N$ is not abelian).}\\

As we have said at the beginning, $x$ is a $p$-element and $y$ is a $q$-element. Suppose that $N/{\rm \textbf{Z}}(N)$ is a $p$-group (the reasoning is analogous if we suppose that it is a $q$-group). Hence we can write $N=P\times A$ where $A\leq {\rm \textbf{Z}}(N)$ and $A$ is a $p'$-group. If $p\neq q$, it follows that $x\in P$ and $y\in A$, which leads to a contradiction with Step 1. Thus, $p=q$ and $x,y \in P$, so the theorem is proved.\\

{\it Step 4. We can suppose that $N$ is not a $\pi$-group.}\\

Let us see that if $N$ is a $\pi$-group, then $N$ is a quasi-Frobenius group with abelian kernel and complement or $N=P\times A$ with $A\leq {\rm \textbf{Z}}(N)$ and $A$ a $p'$-group. Assume that $N$ is a $\pi$-group. As $N$ is non-abelian by Step 3, there exists  a conjugacy class $z^{N}$such that $|z^{N}|\neq 1$, Since $|z^{N}|$ divides $|z^{G}|$, then either $(|z^{N}|, |x^{G}|)=1$ or $(|z^{N}|, |y^{G}|)=1$. Thus, $|z^{N}|$ either is a $\pi_{x}$-number or a $\pi_{y}$-number. If $\Gamma(N)$ is disconnected, we know by Theorem 2 of \cite{BerHerMann} that $N$ is quasi-Frobenius group with abelian kernel and complement. Moreover, $\Gamma(N)$ cannot be empty since by Step 3, we can assume that $N$ is not abelian. Consequently, we can assume that $\Gamma(N)$ is connected and this forces to  either $|x^{N}|=1$ or $|y^{N}|=1$. Suppose for instance that $|x^{N}|=1$, that is, $x\in {\rm\textbf{Z}}(N)$. By Step 3 we can take $w$ an $s$-element of $N\setminus {\rm\textbf{Z}}(N)$ with $s\neq p$. Observe that $|w^{N}|$ must be a $\pi_{y}$-number, so $w^{G}$ is connected to $y^{G}$ in $\Gamma_{G}(N)$. Since $x$ and $w$ have coprime orders and $x\in {\rm\textbf{Z}}(N)$ we have that $|w^{G}|$ and $|x^{G}|$ both divide $|(wx)^{G}|$. As a consequence, we have a contradiction because $|(wx)^{G}|$ has primes in $\pi_{x}$ and $\pi_{y}$. Then we can suppose that $N$ is not a $\pi$-group.\\

{\it Step 5. Conclusion in case $p\neq q$}.\\

Let $z$ be a $\pi'$-element of $K\cap N$ and let us prove that $|z^{G}|$ is a $\pi'$-number. Suppose that $s\in \pi$ is a prime divisor of $|z^{G}|$. We can assume for instance that $s\in \pi_{y}$, otherwise we proceed analogously. Since $|z^{G}|$ divides $|(zx)^{G}|$ we obtain that $s$ divides $|(zx)^{G}|$. On the other hand, we know by the proof of Step 2 that $q\in \pi_{x}$. Therefore, $|(zx)^{G}|$ is divisible by primes in $\pi_{x}$ and $\pi_{y}$, a contradiction. Consequently, $s\not \in \pi $ and $|z^{G}|$ is a $\pi'$-number, as wanted.\\

Let $M$ be the subgroup generated by all $\pi'$-elements of $K\cap N$. Note that $M \neq 1$, otherwise $K\cap N$ would be a $\pi$-group and, since $|N:K\cap N|=|KN:K|$ divides $|G:K|$, which is a $\pi$-number too, then $N$ would be a $\pi$-group, a contradiction with Step 2. Let us prove that $M\unlhd G$. Let $\alpha$ be a generator of $M$, so $|\alpha^{G}|$ is $\pi'$-number. Since $(|G:K|, |\alpha^{G}|)=1$ we have $G=K{\rm\textbf{C}}_{G}(\alpha)$ and hence, $\alpha^{G}=\alpha^{K}\subseteq K\cap N$. Therefore $\alpha^{G}\subseteq M$, as wanted. \\

Let $D=\langle x^{G}, y^{G}\rangle$. Notice that $D\unlhd G$ and $D\subseteq N$. Let $\alpha$ be a generator of $M$. As we have proved that $|\alpha^{G}|$ is $\pi'$-number, then $(|\alpha^{G}|, |x^{G}|)=1$, so $G={\rm\textbf{C}}_{G}(x){\rm\textbf{C}}_{G}(\alpha)$. Thus, $x^{G}=x^{{\rm\textbf{C}}_{G}(\alpha)}\subseteq {\rm\textbf{C}}_{G}(\alpha)$ because $\alpha \in K$. The same happens for $y$, that is, $y^{G}\subseteq {\rm\textbf{C}}_{G}(\alpha)$, so we conclude that $[M,D]=1$.\\

We define $L=MD$ and we distinguish two cases. Assume first that $L< N$. Note that $x, y\in L\unlhd G$ and $L$ trivially satisfies the hypotheses of the theorem. By applying induction to $L$ we have in particular $L={\rm \textbf{O}}_{\pi}(L)\times {\rm \textbf{O}}_{\pi'}(L)$. Observe that the fact that $M\neq 1$ implies that ${\rm \textbf{O}}_{\pi'}(L)> 1$. Now, by definition of $M$, we have that $|K\cap N:M|$ is a $\pi$-number. As $|N:K\cap N|$ is a $\pi$-number, it follows that $|N:{\rm \textbf{O}}_{\pi'}(L)|$ is a $\pi$-number too. Then ${\rm \textbf{O}}_{\pi'}(L)={\rm \textbf{O}}_{\pi'}(N)$ is a Hall $\pi'$-subgroup of $N$. We can apply Lemma 1 so as to conclude that $N={\rm \textbf{O}}_{\pi}(N)\times {\rm \textbf{O}}_{\pi'}(N)$ with $x,y \in {\rm \textbf{O}}_{\pi}(N)$. Since ${\rm \textbf{O}}_{\pi'}(N)> 1$, we apply the inductive hypotheses to ${\rm \textbf{O}}_{\pi}(N)< N$ and we deduce that ${\rm \textbf{O}}_{\pi}(N)$ is a quasi-Frobenius group with abelian kernel and complement or ${\rm \textbf{O}}_{\pi}(N)=P\times A$ with $A\leq {\rm\textbf{Z}}(N)$ and $P$ is a $p$-group, so the theorem is finished. \\

From now on, we assume that $L=N$ and let us see that ${\rm \textbf{Z}}(N)=1$. Otherwise, we take $\overline{N}=N/{\rm\textbf{Z}}(N)$ and $\overline{G}=G/{\rm\textbf{Z}}(N)$. If $|\overline{x}^{\overline{G}}|=1$, then $[\overline{x},\overline{y}]=1$, and thus $[x,y] \in {\rm \textbf{Z}}(N)$. Since $(o(x), o(y))=1$, it is easy to prove that $[x,y]=1$, a contradiction. Analogously, we have  $|\overline{y}^{\overline{G}}|\neq 1$. Consequently, $\overline{N}$ satisfies the assumptions of the theorem. By induction, we have $\overline{N}={\rm \textbf{O}}_{\pi'}(\overline{N})\times {\rm \textbf{O}}_{\pi}(\overline{N})$ with $\overline{x}, \overline{y} \in {\rm \textbf{O}}_{\pi}(\overline{N})$ and ${\rm \textbf{O}}_{\pi}(\overline{N})$ is either a quasi-Frobenius group with abelian kernel and complement or $\overline{N}=\overline{P}\times \overline{A}$ with $\overline{A}\leqslant {\rm \textbf{Z}}(\overline{N})$ and $\overline{P}$ a $p$-group. In the latter case, $[\overline{y}, \overline{x}]=1$ which leads to a contradiction as we have seen before. So we are in the former case. It follows that $N={\rm \textbf{O}}_{\pi'}(N)\times {\rm \textbf{O}}_{\pi}(N)$ with $x,y \in {\rm \textbf{O}}_{\pi}(N)$ and by applying induction to ${\rm \textbf{O}}_{\pi}(N)<N$, we have the result. Therefore, ${\rm\textbf{Z}}(N)=1$. On the other hand, we have proved that $[M,D]=1$. Thus $M\cap D\subseteq {\rm \textbf{Z}}(N)=1$ and $N=M\times D$ with $x,y \in D$. Since $M\neq 1$, we can apply induction to $D$ and get $D={\rm \textbf{O}}_{\pi'}(D)\times {\rm \textbf{O}}_{\pi}(D)$ with $x,y \in {\rm \textbf{O}}_{\pi}(D)$ and ${\rm \textbf{O}}_{\pi}(D)$ is a Frobenius group with abelian kernel and complement (notice that ${\rm \textbf{Z}}({\rm \textbf{O}}_{\pi}(D))=1$ because ${\rm \textbf{Z}}(N)=1$). The $p$-group case cannot occur because $x$ and $y$ do not commute. Notice that if $M$ is $\pi'$-group then the theorem is proved. Assume then that $M$ is not a $\pi'$-group and we will obtain a contradiction. Let $s\in \pi$ such that $s$ divides $|M|$. We can assume that $s\in \pi_{x}$ (we proceed analogously if $s\in \pi_{y}$). Suppose that there exists an $s'$-element $z \in M$ such that $|z^{M}|$ is divisible by $s$. Since $N$ is the direct product of $M$ and $D$, we have that $(zy)^{N}=z^{N}y^{N}$ is a non-trivial class of $N$ whose size is divisible by $s$ and by some prime of $|y^{N}|\neq 1$. This is not possible because $|(zy)^{G}|$ would have primes in $\pi_{x}$ and $\pi_{y}$. Thus, the class size of every $s'$-element of $M$ is a $s'$-number. By Lemma 2, $M=M_{1}\times S$ with $S\in$ Syl$_{s}(M)$. In this case, ${\rm \textbf{Z}}(S)\subseteq{\rm \textbf{Z}}(N)=1$, a contradiction. \\

{\it Step 6. Conclusion in case $p=q$.}\\

Let $K={\rm \textbf{C}}_{G}(x)\cap {\rm \textbf{C}}_{G}(y)$ as in Step 2. Let $z$ be a $p'$-element of $K\cap N$ and let us prove that $|z^{G}|$ is a $\pi'$-number. Suppose that $s\in \pi$ is a prime divisor of $|z^{G}|$. We can assume for instance that $s\in \pi_{y}$, otherwise we proceed analogously. Since $|z^{G}|$ divides $|(zx)^{G}|$ we obtain that $s$ divides $|(zx)^{G}|$. On the other hand, we know by the proof of Step 2 that $q\in \pi_{x}$. Therefore, $|(zx)^{G}|$ is divisible by primes in $\pi_{x}$ and $\pi_{y}$, a contradiction. Consequently, $s\not \in \pi $ and $|z^{G}|$ is a $\pi'$-number, as wanted.\\

Let $T$ be the subgroup generated by all $p'$-elements of $K\cap N$. We have that $T\neq 1$ because otherwise $K\cap N$ would be a $\pi$-group and this implies that $N$ is a $\pi$-group as in Step 5, a contradiction. Let us prove that $T\unlhd G$. If $\alpha$ is a generator of $T$, we know that $|\alpha^{G}|$ is $\pi'$-number. Then $(|G:K|, |\alpha^{G}|)=1$, so we have $G=K{\rm\textbf{C}}_{G}(\alpha)$ and $\alpha^{G}=\alpha^{K}\subseteq K\cap N$. Therefore, $\alpha^{G}\subseteq T$ as wanted. \\

Since the class size of  every $p'$-element of $T$ is a $p'$-number then, by Lemma 2, $T={\rm \textbf{O}}_{p}(T)\times {\rm \textbf{O}}_{p'}(T)$. However, by definition of $T$, we have ${\rm \textbf{O}}_{p}(T)=1$, or equivalently $M={\rm \textbf{O}}_{p'}(T)$. Now, notice that if $s\in \pi$ and $s\neq p$, then the class size of every element of $T$ is an $s'$-number so, it is well known that $T$ has a central Sylow $s$-subgroup and we can write $T={\rm \textbf{O}}_{\pi}(T)\times {\rm \textbf{O}}_{\pi'}(T)$. On the other hand, $|N:T|=|N:K\cap N||K\cap N:T|$ where $|N:K\cap N|=|KN:K|$ is a $\pi$-number and $|K\cap N:T|$ is a power of $p\in \pi$. Therefore ${\rm \textbf{O}}_{\pi'}(T)={\rm \textbf{O}}_{\pi'}(N)$ and ${\rm \textbf{O}}_{\pi'}(N)$ is a Hall $\pi'$-subgroup of $N$. We have proved that the class size of every $p'$-element of $N$ is a $\pi'$-number, so by Lemma 1, we have $N={\rm \textbf{O}}_{\pi'}(N)\times {\rm \textbf{O}}_{\pi}(N)$. We apply induction to ${\rm \textbf{O}}_{\pi}(N)<N$ and the proof is finished. $\Box$ \\

{\it Proof of Corollary B.} It follows immediately by Theorem A. We only have to notice that if $x,y \in {\rm \textbf{O}}_{\pi}(N)$ and $z\in {\rm \textbf{O}}_{\pi'}(N)\setminus {\rm \textbf{Z}}(G)$ it is easy to get that there is a path connecting $x^{G}$ and $y^{G}$ because $(xz)^{G}$ is connected to $x^{G}$ and $(yz)^{G}$ which is connected to $y^{G}$. This contradicts the hypotheses of the Theorem. Thus, ${\rm \textbf{O}}_{\pi'}(N)={\rm \textbf{Z}}(G)$. By the same argument $A\leq {\rm \textbf{Z}}(G)$. $\Box$\\

{\it Proof of Corollary C.} It trivially  follows  by Theorem A. $\Box$\\

We give an example showing that the converse of Theorem A is not true.\\

{\bf Example 1.} We take the Special Linear group $H={\rm SL}(2,5)$ which is a group of order 120 that acts Frobeniusly on $K=\mathbb{Z}_{11} \times \mathbb{Z}_{11}$. Let $P\in$ Syl$_{5}(H)$ and we consider \textbf{N}$_H(P)$. Then, we define $N: = KP$, which is trivially a normal subgroup of $G: =K{\rm \textbf{N}}_H(P)$. We have that the set of the $G$-conjugacy class sizes of $N$ is $\lbrace 1, 20, 242\rbrace$. Thew graph $\Gamma_G(N)$ is form by exactly two vertices joint by an edge. Obviously, $N$ is a Frobenius group with abelian kernel and complement  and there do not exist two non-central $G$-classes in $N$ such that any non-central $G$-class of $N$ has size coprime with one of both.\\

Let us look at two examples illustrating each of the cases of Theorem A.\\

{\bf Example 2}. We take the following groups from the library $SmallGroups$ of GAP (\cite{gap}). Let $G_{1}={\rm Id}(324, 8)$ and $G_{2}={\rm Id}(168, 44)$ that have the normal subgroups that we show now. In fact, $G_{2}$ is the Semi Linear Affine group of order 168. The abelian $3$-subgroup $P=\mathbb{Z}_{3}\times \mathbb{Z}_{3}$ and $A=\mathbb{Z}_{2}\times \mathbb{Z}_{2}\times \mathbb{Z}_{2}$, respectively. It follows that the set of conjugacy class sizes of $P$ is $\lbrace 1, 2, 3, 3\rbrace$ and the set of conjugacy class sizes of $A$ is $\lbrace 1, 7\rbrace$. We construct $N=P\times A$ and $G=G_{1}\times G_{2}$. We have that $N$ is a normal subgroup of $G$ and the set of $G$-conjugacy class sizes of $N$ is $\lbrace 1, 2, 3, 7, 14, 21\rbrace$. So $d(\Gamma_{G}(N))=3$ and $N$ is the direct product of a $3$-group and $A\leq {\rm \textbf{Z}}(N)$. Notice that in this example we have ${\rm \textbf{O}}_{\pi'}(N)=1$ and $\pi=\lbrace 2,3, 7 \rbrace$.\\

{\bf  Example 3}. The quasi-Frobenius case in Theorem A is the natural extension of the ordinary case. It is enough to consider any group $G$ and $N=G$ such that $\Gamma(G)=\Gamma_G(N)$ has two connected components. By the main theorem of \cite{BerHerMann}, we know that $G$ is a quasi-Frobenius group with abelian kernel and complement. \\

\noindent {\bf Acknowledgements}

\bigskip

The results in this paper are part of the third author's Ph.D. thesis at the University Jaume I of Castell\'on.
The research of the first and second authors is   supported by the Valencian Government, Proyecto
PROMETEOII/2015/011. The first and the third authors are also partially supported by Universitat Jaume I, grant P11B\-2015-77.


\begin{thebibliography}{99}

\bibitem{piseparable} A. Beltr\'{a}n and M.J. Felipe, Prime powers as conjugacy class lengths of $\pi$-elements. Bull. Austral. Math. Soc. {\bf 69} (2004), 317-325.

\bibitem{Nuestro1} A. Beltr\'{a}n, M.J. Felipe and C. Melchor, Graphs associated to conjugacy classes of normal subgroups in finite groups. J. Algebra, {\bf 443} (2015), 335-348.

\bibitem{BerHerMann} E.A. Bertram, M. Herzog and A. Mann, On a graph related to conjugacy classes of groups. Bull. London Math. Soc., {\bf 22} (6) (1990), 569-575.

\bibitem{Camina} A.R. Camina, Arithmetical conditions on the conjugacy class numbers of a finite group. J.London Math. Soc. (2) {\bf 5} (1972), 127-132.

\bibitem{Dolfi} S. Dolfi, Arithmetical conditions of the length of the conjugacy classes in finite groups. J. Algebra, {\textbf{174}}, (3) (1995), 753-771.

\bibitem{DolfiCasolo} C. Casolo, S. Dolfi, The diameter of a conjugacy class graph of finite groups.
Bull. London Math. Soc. {\bf 28} (1996), 141-148.

\bibitem{gap} The GAP Group, GAP - Groups, Algorithms and
Programming, Vers. 4.4.12 (2008). http://www.gap-system.org

\bibitem{Kazarin} L.S. Kazarin, On groups with isolated conjugacy classes, Izv.
Vyssh. Uchebn. Zaved. Mat. (7) (1981), 40-45.








\end{thebibliography}
\end{document}